\documentclass[journal, twoside,web]{ieeecolor}
\usepackage{lcsys}
\usepackage{cite}
\usepackage{amsmath}
\usepackage{amstext}
\usepackage{amssymb}
\usepackage{mathtools}
\usepackage{graphicx}
\usepackage{textcomp}
\usepackage{subcaption}
\usepackage{multirow}
\usepackage{booktabs}
\usepackage{algorithm}
\usepackage{algpseudocode}
\usepackage{float}
\usepackage{epsfig}
\usepackage{cite}
\usepackage{tikz}
\usetikzlibrary{fit}
\usetikzlibrary{positioning}
\usetikzlibrary{arrows.meta}
\usepackage{setspace}

\usepackage{fancyhdr}
\newtheorem{theorem}{Theorem}
\newtheorem{assumption}[theorem]{Assumption}

\newtheorem{rem}[theorem]{Remark}
\newtheorem{lemma}[theorem]{Lemma}

\algdef{SE}[PARFOR]{ParFor}{EndParFor}[1]{\algorithmicfor\ #1 \textbf{do in parallel}}{\algorithmicend\ \algorithmicfor}
\newcounter{unnumber}

\newcommand{\R}{\mathbb{R}}

\DeclareMathOperator*\dist{dist}

\def\BibTeX{{\rm B\kern-.05em{\sc i\kern-.025em b}\kern-.08em
    T\kern-.1667em\lower.7ex\hbox{E}\kern-.125emX}}
\begin{document}
\title{Federated Incremental Subgradient Method for Convex Bilevel Optimization Problems}

\author{Sudkobfa Boontawee, Mootta Prangprakhon, and Nimit Nimana
\thanks{This research has received
funding support from the NSRF via the Program Management Unit for Human Resources \& Institutional Development, Research and Innovation [grant number B04G650012]. S. Boontawee is supported by the Development and Promotion of Science and Technology Talents Project (DPST).}
\thanks{The authors are with the Department of Mathematics, Faculty of Science, Khon Kaen University, Khon Kaen 40002, Thailand (e-mail: sudkobfa.b@kkumail.com, mootta\_prangprakhon@hotmail.com, nimitni@kku.ac.th). }
}

\maketitle
\thispagestyle{empty} 

{\color{blue}
\begin{abstract}
In this letter, we consider a bilevel optimization problem in which the outer-level objective function is strongly convex, whereas the inner-level problem consists of a finite sum of convex functions. Bilevel optimization problems arise in situations where the inner-level problem does not have a unique solution.
This has led to the idea of introducing an outer-level objective function to select a solution with the specific desired properties. 
We propose an iterative method that combines an incremental algorithm with a broadcast algorithm, both based on the principles of federated learning. Under appropriate assumptions, we establish the convergence results of the proposed algorithm. To demonstrate its performance, we present two numerical examples related to binary classification and a location problem.
\end{abstract}
}

\section{Introduction}
\label{sec:introduction}

\IEEEPARstart{F}{ederated}  learning (FL) \cite{KMR15, KMRR16, KM16, KM21} is an increasingly popular technique for distributed optimization in machine learning, where multiple clients collaboratively train a model while keeping their data decentralized and without sharing local data with the server or other clients. 
Specifically, each client trains a model using its local data and transmits the updated model to a central server. The central server then collects and aggregates the updates from all local clients to update the global model. This FL scheme supports data stored locally on decentralized clients, ensuring the privacy and security of local data while also reducing the computing and storage burden on the server. 

Due to its numerous advantages, FL is widely employed across various domains. For instance, Sarma \emph{et al.} \cite{SHS21} and Li \emph{et al.} \cite{LMXRHZ19} utilized it in conjunction with deep neural networks for MRI image segmentation. Park \emph{et al.} \cite{PKKKY21} and Ma \emph{et al.} \cite{MZLYLW22} focused on its application in disease diagnosis. Whitmore \emph{et al.} \cite{WMYL25} and Yuxin and Honglin \cite{YH23} applied it to risk control across financial institutions. 
Additionally, Farooq \emph{et al.} \cite{FNMLRSM24}, Guo \emph{et al.} \cite{GZD20}, and Mistry \emph{et al.} \cite{MMSASC23} applied FL in educational contexts, where preserving student data privacy is essential.
For more comprehensive details on FL and its applications, see Banabilah \emph{et al.} \cite{BAAMJ22}, Bharati \emph{et al.} \cite{BMPP22}, and Wen \emph{et al.} \cite{WZL23}, and the references therein.

\vspace{-0.25cm}
\subsection{Problem and Existing Methods}
Typically, FL framework is often framed as minimizing a global loss that aggregates local loss across clients. However, in many practical situations, this can lead to multiple global optimal solutions, and not all of them are equally desirable. One common approach is to introduce an additional criterion that selects the most desirable solution over the set of all global optimal solutions. This idea is the core concept behind the so-called \emph{bilevel optimization problem} (BOP),  which involves two hierarchical levels, namely an outer-level problem and an inner-level problem.

In this letter, we consider a class of BOPs with a structured inner-level problem in the FL setting, formulated as follows: 
\vspace{-0.1cm}
\begin{equation}
\tag{$P^H_F$}
\begin{aligned}\label{problem}
    &\textrm{minimize} \quad && H(x) \\
    &\textrm{subject to} \quad && x \in \textrm{argmin}_{y\in X}F(y),
\end{aligned}
\end{equation}

\noindent
where \(H:\R^n \to \R\) is the outer-level objective function,   $F:\R^n\to\R$ is the inner-level objective function defined as:
\vspace{-0.1cm}
$$F(y):=\sum_{i=1}^{S}  F_i(y),  \text{ where }
F_i(y) := \sum_{j=1}^{I_i} f^i_{j_i}(y),$$ 
\vspace{-0.35cm}

\noindent and $X\subset\R^n$ is the constraint. Here, $S$ represents the total number of clients and each function $F_i$ corresponds to the local loss associated with the client $i$. 
Moreover, each local client loss $F_i$ is defined as the sum of losses over its local data points, with \(f^i_{j_i}\) denoting the loss associated with the \(j\)-th data point on the client \(i\), and $I_i$ denotes the number of local data points associated with the client $i$. 

For simplicity of notation, we let $[s] := \{ 1,2,\ldots,s\}.$ We first make the following assumption on Problem (\ref{problem}):

\begin{assumption}\label{ass1}  
\begin{itemize}    
    \item[(i)] \hspace{-0.1cm}The function $H : \mathbb{R}^n \rightarrow \mathbb{R}$ is $\mu_H$-strongly convex with $\mu_H>0$.
    
    \item[(ii)] \hspace{-0.1cm}The functions $f_{j_i}^i : \mathbb{R}^n \rightarrow \mathbb{R}$, $i\in [S],j_i\in [I_i]$, are convex. 

     \item[(iii)] \hspace{-0.1cm}The set $X \subset \mathbb{R}^n$ is  nonempty, compact, and convex.
\end{itemize}
\end{assumption}

We denote by $x^*_{H}$ the unique solution to Problem (\ref{problem}) and $m:=\sum_{i=1}^S I_i$ the total number of inner-level functions.

{\color{black}
In the literature, several approaches have been proposed to address various classes of Problem (\ref{problem}) across different contexts. 
A standard and widely used approach involves solving the following regularized problem: 
 \begin{equation}
\tag{$P_{\lambda}$}
\begin{aligned}\label{problemlambda}
    &\textrm{minimize} \quad && F(x) + \lambda H(x) \\
    &\textrm{subject to} \quad && x \in X.
\end{aligned}
\end{equation}
An example of this regularization technique is the well-known {\it Tikhonov regularization} \cite{TA77}. 
Under some restrictive conditions, it has been shown that there is a sufficiently small $\lambda^*>0$ such that the optimal solution to Problem (\ref{problemlambda}) is also the optimal solution of Problem (\ref{problem}). Notwithstanding, in practice, the exact value $\lambda^*$ is not known. 

To overcome this, the {\it iterative regularization} technique has been proposed, in which the strongly convex case of Problem (\ref{problemlambda}) is studied by replacing the parameter $\lambda$ with a sequence of positive regularization parameters $\{\lambda_k\}_{k=1}^\infty$. 
We denote Problem (\ref{problemlambda}) with $\lambda=\lambda_k$ as Problem $(P_{\lambda_k})$ for all $k\geq1$. By denoting the unique minimizer of Problem $(P_{\lambda_k})$ by $x_{\lambda_k}^*$, one can ensure that the sequence $\{ x^* _{\lambda_k}\}_{k=1}^\infty$ converges to $x^* _H$ (see Lemma \ref{lemma1} below). 
This iterative regularization technique originates by Solodov \cite{S07}, who proposed the {\it explicit descent method}. 
 Building on this approach,  Amini and Yousefian \cite{AY19} presented an {\it iterative regularized incremental projected subgradient method} (IR-IG), in which they incorporated the ideas of the incremental projected subgradient method \cite{NB01} with iterative regularization. IR-IG is designed for solving Problem (\ref{problem}) with $S=1$ (i.e., $I_1=m$) and its general step is defined as follows: set $x_{k,1} = x_k$, then compute 
 \vspace{-0.1cm}
 \begin{align}\label{irig} x_{k,j+1} = P_{X}\left[ x_{k,j} - \gamma_k g_{k,j} - \frac{\gamma_k \lambda_k}{m} \mathcal{H}_{k,j}\right],
 \end{align}
for all $j\in [m]$ and update  $x_{k+1}$ as 
$x_{k+1} = x_{k,m+1},$
where $ g_{k,j}$ is a subgradient of $f_{j}^1 $ at $x_{k,j}$,  $\mathcal{H}_{k,j}$ is a subgradient of $H$ at $x_{k,j},$ and  $\gamma_k,\lambda_k$ are positive stepsizes. Note that, in many situations, direct communication between clients as in IR-IG may not be suitable due to data privacy concerns. This limitation  
leads to the consideration 
of FL settings. Recently, Ebrahimi et al.\cite{EQC25} proposed a universal regularized framework for smooth BOPs in the FL setting. In their approach, each client solves Problem (\ref{problemlambda})
using its corresponding inner-level and outer-level objective functions. After all clients complete their local updates, the central server aggregates the results and updates the global model accordingly.

\vspace{-0.25cm}
\subsection{Contributions}

The main contribution of this letter is to present the {\it federated incremental subgradient method} (FISM), designed to solve Problem (\ref{problem}). This method is a combination of IR-IG  and the concept of FL framework.  
To be specific, the method performs systematic computations similar to those in FL context, and employs an incremental projected subgradient scheme to update each local client’s model during their local training steps. Note that IR-IG requires the subgradient of the outer-level objective function $H$ to be computed at every local update, whereas the method allows the subgradient of $H$ to be performed only once per global update step. 
This less frequent computation may reduce computational demands and may also decrease communication overhead.
Under certain assumptions, we establish both the convergence of the method and the convergence rate of the inner-level function value of the average iterate toward the optimal value of Problem (\ref{problem}). 
Moreover, to demonstrate its effectiveness and performance, we present two numerical examples and compare the proposed method’s performance with that of IR-IG.

{\it Notations:} Throughout this letter,  $\R^n$ denotes a Euclidean space equipped with the inner product $\langle \cdot, \cdot \rangle$ and norm $\|\cdot\|$. 
Let $C\subset\R^n$ and $x\in\R^n$ be given. We denote the {\it distance} from $x$ to $C$ by $\textrm{dist}(x,C):=\textrm{inf}_{c\in C}\|x-c\|$. If $C$ is nonempty, closed and convex, then there exists a unique $y\in C$ such that $\|x-y\|=\dist(x,C)$. We call $y$ the {\it metric projection} of $x$ onto $C$ and denote it by $P_Cx$. Let  $f:\R^n\to\R$ be a function.
The function $f:\R^n\to\R$  is  {\it convex} 
if $f(\alpha x+ (1-\alpha)y)\le \alpha f(x)+(1-\alpha)f(y)$ for all $x,y\in\R^n$ and $\alpha\in[0,1]$; and
{\it $\rho$-strongly convex}, where $\rho>0$, if 
$f(\alpha x+ (1-\alpha)y)\le \alpha f(x)+(1-\alpha)f(y)-\frac{\rho}{2}\alpha(1-\alpha)\|x-y\|^2$.
If the function $f$ is convex, we define a {\it subgradient} of $f$ at $x\in\R^n$ to be a vector $g\in\R^n$ in which 
$\langle g,y-x\rangle \le f(y)-f(x)$ for all $y\in\R^n$. The set of all subgradients of  $f$ at $x$ is called the {\it subdifferential} of $f$ at $x$, and is noted by
$\partial f(x)$.
If the function $f$ is $\rho$-strongly convex, then for all $x,y\in\R^n$, it holds that $f(y)\ge f(x)+\langle g,y-x\rangle + \frac{\rho}{2}\|y-x\|^2$, where $g\in \partial f(x)$.

\section{Federated Incremental Subgradient Method}\label{section-algorithm}

Next, we present an iterative algorithm to solve Problem (\ref{problem}) and investigate its convergence. 

\vspace{-0.2cm}

\begin{algorithm}[H]
\caption{{\small FISM: Federated Incremental Subgradient Method}}
\label{alg}
\begin{algorithmic}[1]
    \State \textbf{Input:} Starting point $x_1 \in \mathbb{R}^n$, 
    positive stepsizes $\{ \gamma_k \}_{k=1} ^ \infty$ and $\{\lambda_k \}_{k=1}^\infty $.
    \For{$k = 1, 2, \dots$},
        \State Pick $\mathcal{H}_k \in \partial H(x_k)$.
        \State Send $x_k$ and $\mathcal{H}_k$ to clients.
        \For{all clients $i \in [S]$, \textbf{in parallel}}
            \State $ x_{k,1}^{i} = x_k$. 
            \For{all local users $j_i \in [I_i]$}
            \State Pick $g_{k,j_i}^i \in \partial f_{j_i}^i (x_{k,j_i}^i)$.
            \State {\small $ x_{k,j_i+1}^{i} = P_X \left[ x_{k,j_i}^{i} - \gamma_k g_{k,j_i}^i -  \frac{\gamma_k \lambda_k}{m}  \mathcal{H}_k \right]$. }
            \EndFor
            \State $x_{k}^{i} = x_{k,I_i + 1} ^{i}$.
        \EndFor
        \State $ x_{k+1} = \frac{1}{S} \sum _{i=1}^{S} x_k ^{i}$.
    \EndFor
\end{algorithmic}
\end{algorithm}
\vspace{-0.4cm}
The workflow of FISM is presented in Fig. \ref{framework-algor}.
\vspace{-0.4cm}

\begin{figure}[H]
\resizebox*{9.2cm}{!}{\includegraphics{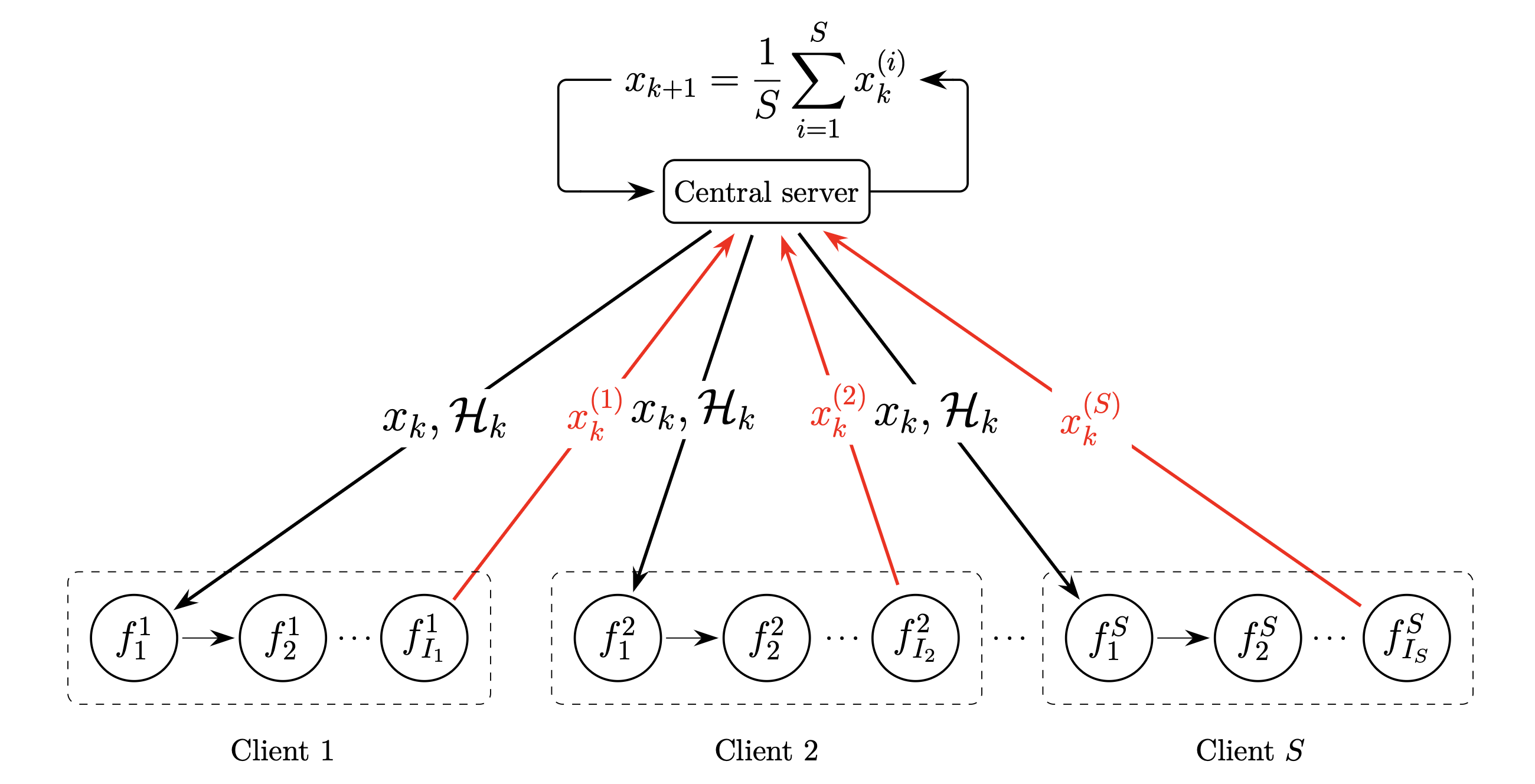}}
\vspace{-0.3cm}
\caption{Workflow of FISM}
\label{framework-algor}
\end{figure}
\vspace{-0.3cm}

As shown in FIg \ref{framework-algor}, 
the central server firstly distributes the parameters $x_k$ and $\mathcal{H}_k$ to all clients, then each client performs local training on their private data,  and subsequently sends their updated model parameters back to the central server.

Next, we provide some remarks regarding FISM.

\begin{rem}\label{remark3.1}
Notice that when $S=1$, FISM is closely related to IR-IG. However, a key distinction between them lies in the term $\mathcal{H}_k$, which appears in FISM. This term allows the central server to preserve the privacy of the outer-level objective function $H$ by sharing only the subgradient $\mathcal{H}_k$ with all clients, rather than revealing the full explicit form of $H$. 
As a result, the central server can protect the privacy of its outer-level objective function. This feature further distinguishes FISM from the existing framework proposed by Ebrahimi et al. \cite{EQC25}.
In addition, since each local user uses only the subgradient $\mathcal{H}_k$ provided by the central server, rather than computing new subgradients locally, FISM reduces the computational cost per local step. 
Specifically, in each iteration where all data points are fully visited, 
 IR-IG requires computing subgradients for $2m$ times, whereas FISM reduces this to $m+1$ times.
\end{rem}

\begin{rem}\label{remark3.3} 
    Another aspect,
    FISM and IR-IG also differ in their total computational running time per iteration.
    Note that since FISM is based on the FL setting, it involves communication costs between the central server and all clients, which can affect computational running times.
    Indeed, let \( t_{j} \) and \( s_{i,j_i} \)
    denote the computational times required for each local update corresponding to the local data points in IR-IG and FISM, which correspond to formula (\ref{irig}) and  line 9 of Algorithm \ref{alg}, respectively.  
    We let the communication time between the central server and client $i$  be $\varepsilon_i$. In this case, the times required per iteration for IR-IG and FISM satisfy
    $T_{\text{IR-IG}}  = \sum_{j = 1}^{m} t_{j}$ and $T_{\text{FISM}} \leq \textrm{max}_{i \in [S] }\sum_{j_i = 1}^{I_i} s_{i,j_i} +\textrm{max}_{i \in [S] }\varepsilon_i,$
    where $T_{\text{IR-IG}}$ and $T_{\text{FISM}}$ represent the times required per iteration for IR-IG and FISM, respectively. Moreover, the relation between $s_{i,j_i}$ and $t_j$ satisfies $s_{i,j_i} \leq t_j$, which implies that 
\vspace{-.1cm}
    $$T_{\text{FISM}} \leq T_{\text{IR-IG}}+\textrm{max}_{i \in [S] }\varepsilon_i.$$
    \vspace{-0.35em}
    
\end{rem}

\vspace{-0.15cm}

}

\vspace{-0.25cm}
Next, we shall discuss some useful properties as follows: by invoking the boundedness of \( X \) and \cite[Proposition 16.20]{BC17}, we get that the subdifferentials of \( H \) and \( f^i_{j_i} \) on \( X \) are bounded, i.e., there exist constants \( B_H, B_f > 0 \) such that
$\| g_h \| \leq B_H$ and $\| g_{j_i}^i \| \leq B_f$, where $g_h\in\partial H(x)$ and $g_{j_i}^i\in\partial f^i_{j_i}(x)$
for all \( i \in [S] \), \( j_i \in [I_i] \), and for all \( x \in X \). Moreover, there exist constants \( L_f, L_H > 0 \) such that
$|f_{j_i}^i(x) - f_{j_i}^i(y)| \leq L_f \| x - y \|$ and $|H(x) - H(y)| \leq L_H \| x - y \|$,
for all \( i \in [S] \), \( j_i \in [I_i] \), and for all \( x, y \in X \).  Furthurmore, we know that $H$ is also bounded on $X$, i.e., there is a constant $M_H > 0$ such that $| H(x) | < M_H $ for all $x \in X$.
Thus, throughout this letter, we set $C_f := \textrm{max}\{ B_f,L_f \}$ and $C_H := \textrm{max}\{B_H,L_H,M_H \}.$

The following lemma provides a key tool in proving the convergence analysis of FISM.

\begin{lemma}\cite[Proposition 1]{YNS17} \label{lemma1}
Let Assumption \ref{ass1} hold, and let $\{ \gamma_k \}_{k=1}^\infty$ be a positive sequence. Let $x_{\lambda_k}^*$ be the minimizer of Problem $(P_{\lambda_k})$  for all $k \geq 1$. Then,
\begin{itemize}
\item[(i)] For every $k \geq 2,$  $\| x_{\lambda_k}^* - x_{\lambda_{k-1}}^* \| \leq \frac{C_H}{\mu_H} \left| 1 - \frac{\lambda_{k-1}}{\lambda_k} \right|.$
\item[(ii)] If $ \lambda_k \rightarrow 0,$ then  $x_{\lambda_{k}}^* \rightarrow x_H ^* $ .
\end{itemize}
\end{lemma}

Next, we present a useful inequality that will play an important role in the sequel.
For the sake of compactness, we denote $V_{k}(x):=\|x_k - x\|^2$ and $V_{k,j_i } ^i (x):=\|x_{k,j_i } ^i - x\|^2$, for all $i \in [S]$, $j_i \in [I_i]$, $x\in X$ and $k\geq1$.

\vspace{0.15cm}
\begin{lemma}\label{lemma5} Let $\{x_k\}_{k=1}^\infty$ be a sequence generated by FISM. Suppose that Assumption \ref{ass1} holds, $0 < \gamma_k \lambda_k\mu_H  \leq 2m$ for all $k\geq1$, and $\{\lambda_k \}_{k=1}^\infty $ is a nonincreasing sequence. Then, for all $k \geq 2$, the following inequality holds:
{\small \begin{align*}
V_{k+1} (x_{\lambda_k}^*)
\leq& \left( 1 - \frac{\gamma_k \lambda_k \mu_H}{2m} \right) V_{k} (x_{\lambda_{k-1}}^*) \\&+  \frac{4m C_H ^2}{\gamma_k \lambda_k \mu_H ^3} \left( 1 - \frac{\lambda_{k-1}}{\lambda_k} \right)^2 + \gamma_k ^2 C,
\end{align*} }
\noindent where $C :=  4m^3 \left( C_f ^2 + \frac{2\lambda_1 ^2 C_H ^2}{m^2} \right) + 2 C_f m^2  \left(C_f + \frac{\lambda_1 C_H }{m}\right)$ and $x_{\lambda_k} ^*$ is the unique minimizer of Problem $(P_{\lambda_k})$.
\end{lemma}
\begin{proof}
    Let $k\geq 2$ be given. By the definition of $x_{k,j_i + 1}^{i}$ and the nonexpansiveness of the metric projection $P_X$, we have, for all $x\in X$, $i \in [S]$, $j_i \in  [I_i]$, that
{\small \begin{align*}
   &V_{k,j_i +1} ^i (x) 
    \\ \leq &  V_{k,j_i } ^i (x) + \gamma_k ^2 \left\| g_{k,j_i}^i - \frac{ \lambda_k}{m} \mathcal{H}_k \right\|^2  - 2\gamma_k \left\langle g_{k,j_i}^i , x_{k,j_i}^{i} - x \right \rangle \nonumber
    \\  &- \frac{2\gamma_k \lambda_k}{m} \left\langle \mathcal{H}_k, x_{k,j_i}^{i} - x \right\rangle \nonumber
    \\ \leq &  V_{k,j_i } ^i (x) + 2\gamma_k ^2 \left( B_f ^2 + \frac{ \lambda_k ^2 B_H ^2}{m^2} \right) -2\gamma_k \left(f_{j_i}^i (x_{k,j_i}^{i}) - f_{j_i}^i (x)\right)  \nonumber
    \\& - \frac{2\gamma_k \lambda_k}{m} \left[ 
    \left\langle \mathcal{H}_k, x_{k,j_i}^{i} - x_k \right\rangle  
    - \left( H(x) - H(x_k) - \frac{\mu_h}{2} V_k(x) \right)
\right],
\end{align*}}
\noindent where the last inequality follows from the boundedness of subgradients, and the characterization of strong convexity.
Thus, by applying Young's inequality to the term involving inner product, adding and subtracting the term $f^i_{j_i}(x_k)$, and rearranging, the above inequality becomes 
{\small \begin{align}\label{eq1}
    &V_{k,j_i + 1}^i (x) \nonumber    \\
    \leq& V_{k,j_i}^i (x) +  2\gamma_k ^2\left( C_f ^2 + \frac{\lambda_k ^2C_H ^2}{m^2} \right)  + \frac{2 \gamma_k^2 \lambda_k^2  C_H ^2}{m^2}  +\| x_{k,{j_i}}^{i} - x_k \|^2 
    \nonumber
    \\& - 2\gamma_k \left( f_{j_i}^i (x_k) + \frac{\lambda_k}{m}H(x_{k}) \right) + 2\gamma_k \left(f_{j_i}^i (x) + \frac{\lambda_k}{m}H(x) \right)
    \nonumber
    \\& - \frac{\gamma_k \lambda_k \mu_h}{m} V_k(x) + 2 \gamma_k \left( f_{j_i}^i (x_k) - f_{j_i}^i (x_{k,j_i}^{i}) \right),
\end{align}}
\noindent for all $x\in X$, $i \in [S]$, $j_i \in [I_i]$. To investigate the last term of (\ref{eq1}), we note that
{\small \begin{align*} 
    \| x_{k,j_i}^{i} - x_k \| 
    \leq \sum_{l = 1}^{j_i -1} \left\| x_{k,l+1}^{i} - x_{k,l}^{i}\right\|
    \leq j_i\gamma_k \left(C_f + \frac{\lambda_k C_H}{m}  \right), 
\end{align*} }
which implies that, for all  $i \in [S]$, $j_i \in [I_i]$,
\begin{align}\label{lastremark}
    \| x_{k,j_i}^{i} - x_k \|^2  \leq 2m^2 \gamma_k ^2 \left( C_f ^2 + \frac{\lambda_k^2C_H ^2}{m^2} \right).
\end{align}
Since $\{\gamma_k \}_{k=1}^\infty $ is nonincreasing, we set
$$A_1:= C_f ^2 + \frac{2\lambda_1 ^2 C_H ^2}{m^2} \quad \text{and} \quad A_2:=C_f + \frac{\lambda_1 C_H}{m}.$$ 
Thus, these relations, together with the Lipschitz continuity of $f_{j_i}^i$ over $X$, imply that  inequality (\ref{eq1}) becomes 
{\small \begin{align*}
    & V_{k,j_i + 1}^i (x) \\
    \leq& V_{k,j_i}^i (x) - \frac{ \gamma_k \lambda_k \mu_h}{m} V_k (x) +  2(m^2+1)\gamma_k ^2 A_1 + 2j_i \gamma_k ^2 C_f A_2
    \nonumber
    \\& - 2\gamma_k \left( f_{j_i}^i (x_{k}) + \frac{\lambda_k}{m}H(x_{k}) \right) + 2\gamma_k \left(f_{j_i}^i (x) + \frac{\lambda_k}{m}H(x) \right),
    \nonumber
\end{align*}  }
\noindent for all $x\in X$, $i \in [S]$, $j_i \in [I_i]$. On each client $i$, summing this inequality over the number of local users  from $j_i=1$ to $I_i$, telescoping the terms, and invoking the definition $V_{k,1}^i (x)=V_{k} (x)$, along with the fact that $I_i\leq m$, $\sum_{j_i=1}^{I_i} j_i=\frac{I_i(I_i+1)}{2}\leq m^2$, we obtain, for all $x\in X$, $i \in [S]$, that
{\small \begin{align*} 
    &V_{k,I_i + 1}^i (x) \\
    \leq& \left( 1 - \frac{I_i \gamma_k \lambda_k \mu_H}{m} \right) V_{k} (x) + 4m^2 I_i\gamma_k ^2 A_1 + 2m^2\gamma_k ^2 C_f A_2
    \\& - 2\gamma_k  \left( F_i(x_{k}) + \frac{\lambda_k I_i}{m}H(x_{k}) -  F_{i}(x) + \frac{I_i \lambda_k}{m}H(x) \right)\nonumber. 
\end{align*} }
\noindent Summing this inequality over the number of clients from $i=1$ to $S$, and  dividing by $S$, we obtain for all $x\in X$ that
{\small \begin{align}\label{rateeq}
    V_{k+1} (x)\leq& \left( 1 - \frac{\gamma_k \lambda_k \mu_H}{m} \right) V_k (x) + \gamma_k ^2\left(4m^3 A_1 + 2 m^2 C_f A_2\right) \nonumber
    \\& + \frac{2\gamma_k}{S} \left( (F + \lambda_kH)(x) - (F +  \lambda_kH)(x_k) \right) ,
\end{align} }
\noindent where the left-hand side follows from Jensen's inequality.
Thus, substituting $x$ by $x_{\lambda_k}^*$, and using the facts  $x_{\lambda_k}^*$ is the unique minimizer of $F+\lambda_kH$ over $X$, and that $x_k\in X$, the above inequality becomes
{\small \begin{align}\label{eq2}
    &V_{k+1} (x_{\lambda_k}^*)  
    \leq \left( 1 - \frac{\gamma_k \lambda_k \mu_H}{m} \right)V_k (x_{\lambda_k}^*) + \gamma_k ^2 C .
\end{align}}
Again, by invoking Young's inequality, we note that
{\small \begin{align*}
     V_k (x_{\lambda_k}^*)
    \leq& \left(1+\frac{\gamma_k \lambda_k \mu_H}{2m}\right) V_k (x_{\lambda_{k-1}}^*) \\ &+ \left(1 + \frac{2m}{\gamma_k \lambda_k \mu_H}\right)\| x_{\lambda_{k-1}}^* - x_{\lambda_k} ^* \|^2, 
\end{align*} }
and by substituting this relation into  (\ref{eq2}), we obtain
\vspace{-0.15cm}
{\small \begin{align*}
&V_{k+1} (x_{\lambda_k}^*) \nonumber\\
     \leq& \overbracket{\left( 1 - \frac{\gamma_k \lambda_k \mu_H}{m} \right)\left(1+\frac{\gamma_k \lambda_k \mu_H}{2m}\right)}^{(a)} V_k (x_{\lambda_{k-1}}^*) \nonumber 
     \\&+ \overbracket{\left( 1 - \frac{\gamma_k \lambda_k \mu_H}{m} \right)\left(1 + \frac{2m}{\gamma_k \lambda_k \mu_H}\right)}^{(b)}\| x_{\lambda_{k-1}}^* - x_{\lambda_k} ^* \|^2 + \gamma_k ^2 C.\nonumber
\end{align*} }

\noindent Note that $(a)\le 1 - \frac{\gamma_k \lambda_k \mu_H}{2m},$ and $(b)\le 1+ \frac{2m}{\gamma_k \lambda_k \mu_H} \le \frac{4m}{\gamma_k \lambda_k \mu_H},$ where the final bound holds  due to $\frac{2m}{\gamma_k \lambda_k \mu_H} \geq 1$.
Hence, by associating these terms and invoking Lemma \ref{lemma1}(i), we obtain the required inequality.
\end{proof}

\vspace{0.2cm}
Next, we adopt the following assumption which is proposed by Amini and Yousefian \cite[Assumption 2]{AY19}.

\begin{assumption}\label{ass2}
The following conditions are assumed:
\begin{itemize}
\item[(i)] $\{ \gamma_k \}_{k=1}^\infty$ and $\{ \lambda_k \}_{k=1}^\infty$ are nonincreasing sequences of positive real numbers with $\gamma_1 \lambda_1 \mu_H \leq 2m$;
\item[(ii)] $\sum_{k=1}^\infty \gamma_k \lambda_k = \infty, \sum_{k=1}^{\infty} \frac{1}{\gamma_{k+1} \lambda_{k+1}} \left( 1 - \frac{\lambda_{k}}{\lambda_{k+1}} \right)^2 < \infty,\\ \sum_{k =1}^\infty \gamma_k ^2 < \infty$;
\item[(iii)] $ \frac{1}{\gamma_{k+1} ^2 \lambda_{k+1} ^2} \left( 1 - \frac{\lambda_{k}}{\lambda_{k+1}} \right)^2 \rightarrow 0,\;  \frac{\gamma_k}{\lambda_k} \rightarrow 0,\; \lambda_k \rightarrow 0. $
\end{itemize}
\end{assumption}
\vspace{0.2cm}

\begin{rem}\label{remarkss}
    An example of choices for the stepsizes $\{ \gamma_k \}_{k=1}^\infty$ 
    and $\{ \lambda_k \}_{k=1}^\infty$ that satisfy Assumption \ref{ass2} is $\gamma_k = \frac{\gamma_1}{k^a}$ and $\lambda_k = \frac{\lambda_1}{k^b}$ for all $k \geq 1$, where $\gamma_1,\lambda_1 > 1, \gamma_1 \lambda_1 \mu_h \leq 2m, a > 0.5, a>b>0 $, and $a+b < 1$.
\end{rem}

\vspace{0.2cm}

We are now in a position to prove the main theorem. Before doing so, we present another useful lemma, which will play a key role in the convergence analysis.
\vspace{.2cm}
\begin{lemma}\cite[pp. 49]{P87}
\label{P87-Lm10}
Let $\{a_k\}_{k=1}^\infty$ be a sequence of nonnegative real numbers, and let $\{\alpha_k\}_{k=1}^\infty$ and $\{\beta_k\}_{k=1}^\infty$ be sequences of nonnegative real numbers such that for all $k\ge1$, 
$a_{k+1}\le (1-\alpha_k)a_k + \beta_k,$ 
where
$\alpha_k\in(0,1]$, $\sum_{k=1}^\infty\alpha_k=\infty$, $\sum_{k=1}^\infty\beta_k<\infty$, and $\frac{\beta_k}{\alpha_k}\to0$. Then, $a_k\to 0$. 
\end{lemma}

\begin{theorem}\label{converge}
Let $\{ x_k \}_{k=1}^\infty$ be a sequence generated by FISM. Suppose that Assumptions \ref{ass1} and \ref{ass2} hold. Then, $\{ x_k \}_{k=1}^\infty$ converges to  $x_H ^*$.
\end{theorem}
\begin{spacing}{1}
\begin{proof} For all $k \geq 1,$ we set 
$\alpha_k = \frac{\gamma_{k+1} \lambda_{k+1} \mu_H}{2m}$
and  $ \beta_k =  \frac{4m C_H ^2}{\gamma_{k+1} \lambda_{k+1} \mu_H ^3} \left( 1 - \frac{\lambda_{k}}{\lambda_{k+1}} \right)^2 + \gamma_k ^2 C,$
and $a_k:=V_{k+1}(x^*_{\lambda_k})$.
By Assumption \ref{ass2}(i) and (ii), we obtain that $\alpha_k \in[0,1]$ and $\sum_{k=1}^{\infty} \alpha_k = \infty$. Moreover, the assumption that $\{ \lambda_k \}_{k=1}^\infty$ is nonincreasing
together with Assumption \ref{ass2}(ii), 
imply that $\sum_{k=1}^\infty \beta_k < \infty.$
Furthermore,
we also have 
 $\gamma_k \lambda_k \leq \gamma_k \lambda_1 \leq \frac{\lambda_1 ^2 \gamma_k}{\lambda_k},$ which provides an upper bound for
 $\{\gamma_k \lambda_k\}_{k=1}^\infty$. By this and Assumption \ref{ass2}(iii), we obtain that $\frac{\beta_k}{\alpha_k} \rightarrow 0.$
Applying Lemma \ref {P87-Lm10},
we obtain  $ \|x_{k+1} - x_{\lambda_{k}} ^*\|\to0$.
 Hence, by invoking Assumption \ref{ass2}(iii), Lemma \ref{lemma1}(ii), and the triangle inequality, we conclude that $ x_{k+1} \rightarrow x_H ^* ,$
which completes the proof.
 \end{proof}   
\end{spacing}

\vspace{-0.025cm}
\begin{rem}
By applying inequality (\ref{lastremark}), 
the obtained result in Theorem \ref{converge}, and the triangle inequality, it follows that the sequence $\{ x_{k,I_i + 1}^{i} \}_{k=1}^\infty $ also converges to the unique solution of Problem (\ref{problem}) for all $i \in [S]$.  This result implies that each client can access the unique solution of Problem (\ref{problem}). Therefore, it aligns well with the concept of the FL, which enables each client to preserve privacy while still contributing to the solution of the system.
\end{rem}

By choosing particular  stepsizes, we obtain the convergence rate of the inner-level function value of the average
iterate to the optimal value of Problem (\ref{problem}) as follows.

\begin{theorem}
Let $\{ x_k \}_{k=1}^\infty$ be a sequence generated by FISM. Suppose that Assumption \ref{ass1} holds and $\epsilon \in (0,0.5)$. 
If the sequences $\{\gamma_k\}_{k=1}^\infty$ and $\{\lambda_k\}_{k=1}^\infty$ are defined by $$\gamma_k = \frac{\gamma_1}{k^a} \quad \text{and}  \quad \lambda_k 
= \frac{\lambda_1}{k^b},$$ 
where $\gamma_1,\lambda_1 > 0 $ with $\gamma_1 \lambda_1 \mu_H \leq 2m$,  $a = 0.5 + 0.5\epsilon, b = 0.5-\epsilon$,
then we have for $K\geq 1$ that
{\small\begin{align*}
    F(\hat{x}_K) - F(x^* _H) \leq \mathcal{O} \left( \frac{1}{K^{0.5 - \epsilon}} \right),
\end{align*}}
\noindent where $\hat{x}_K := \frac{\sum_{k=1}^K \gamma_k x_k}{\sum_{k=1}^K \gamma_k}$.
\end{theorem}

\begin{proof}
 Putting $x:= x^* _H$ in (\ref{rateeq}), we obtain 
{\small \begin{align*}
    V_{k+1} (x^* _H)\leq& \left( 1 - \frac{\gamma_k \lambda_k \mu_H}{m} \right) V_k (x^* _H) + \gamma_k ^2C \nonumber
    \\& + \frac{2\gamma_k}{S} \left( (F + \lambda_kH)(x^* _H) - (F +  \lambda_kH)(x_k) \right)
    \\ \leq& V_k (x^* _H) + \gamma_k ^2 C 
    \\& + \frac{2\gamma_k}{S} \left( (F(x^* _H) - (F(x_k)) + 2 \lambda_k C_H)\right),
\end{align*}}

\noindent where $C$ is given in Lemma \ref{lemma5} and the last inequality follows from the boundedness of $H$.
Summing this inequality from $k = 1$ to $N$, we obtain
\vspace{-.2cm}
{\small \begin{align*}
    \sum_{k = 1}^N \frac{2\gamma_k}{S} \left(F(x_k) - F(x^* _H)\right) \leq \frac{4 C_H}{S} \sum_{k=1}^N \gamma_k \lambda_k + M,
\end{align*}}

\noindent where $M := V_0 (x^* _H) + C \sum_{k=1}^\infty \gamma_k ^2.$ Dividing
this inequality by $\frac{2 \sum_{k=1}^N \gamma_k }{S}$  and applying Jensen's inequality, we have
{\small
\begin{align*}
    F(\hat{x}_K) - F(x^* _H) \leq \frac{4 C_H \sum_{k=1}^N \gamma_k \lambda_k + SM}{2 \sum_{k=1}^N \gamma_k}.
\end{align*}
}

\noindent Investigating the  upper bounds {\small
$$ \frac{1}{\sum_{k=1}^N \gamma_k} 
\leq \mathcal{O} \left(\frac{1}{K^{1-a}}\right),
\textrm{ and } 
  \frac{\sum_{k=1}^N \gamma_k \lambda_k}{\sum_{k=1}^N \gamma_k} \leq \mathcal{O}\left( \frac{1}{K^{b}}\right),$$
}
together with the choices of $a$ and $b$, 
we infer that
\begin{align*}
    F(\hat{x}_K) - F(x^* _H) \leq \mathcal{O} \left( \frac{1}{K^{b}} \right),
\end{align*}
which completes the proof.
\end{proof}

\vspace{-0.2cm}

\section{Numerical Experiments}
\label{section-experiment}
\subsection{Binary Image Classification}
 
\vspace{-0.1cm}

In this experiment, we compare the performance of FISM and IR-IG using the MNIST dataset \cite{LBB98}, focusing on binary classification of handwritten digits $0$ and $1$. The problem can be written as: 
\vspace{-0.1cm}
{\small \begin{equation*}
\begin{aligned}\label{expirement}
    &\textrm{minimize} \quad &&  \|x\|_1 + 0.5\|x\|^2 \\
    &\textrm{subject to} \quad && x \in {\textrm{argmin}} _{y \in X}  {\small \sum_{i=1} ^S \sum_{j_i=1}^{I_i} f_{j_i} ^i(y)},
\end{aligned}
\end{equation*}
}
\vspace{-0.2cm}

\noindent where $S$ is the number of groups (clients), $I_i$ is the number of samples in group $i$, and $f_{j_i}^i := \log(1+\exp{(-b_{j_i}^i \langle a_{j_i}^i,\cdot \rangle)})$ is the logistic loss of function corresponding to $j_i -$th sample $(a_{j_i}^i,b_{j_i}^i) \in \mathbb{R}^{784} \times \{ \pm{1}\}$ in group $i$,\, and $X := [-100,100]^{784}$ is the box.  
In this problem, the lower-level problem corresponds to logistic regression, which may admit multiple solutions, while the outer-level problem selects a single solution that simultaneously provides sparsity and minimal norm.

First, we demonstrate the performance of FISM with different numbers of groups ($S=1,2,4,8$) compared to IR-IG. We focus on the effects of the total number of observations ($m$) and the computational running time (in seconds).

To compare, we choose the possibly best fine-tuned parameters of IR-IG, which are $ \gamma_k := \frac{10}{k^{0.8}}$  and $\lambda_k := \frac{1}{k^{0.1}}$, and use them when evaluating both  FISM and IR-IG.
We set the maximum number of iterations to be $200$. The experiments are conducted using $m$ training samples, with each class containing $m/2$ samples. We show the averaged results over $10$ independent runs in which both the initial points and the ordering of the training data are randomized.
 \vspace{-0.2cm}
\begin{figure}[h!]
    \centering
\includegraphics[width=0.6\linewidth ]{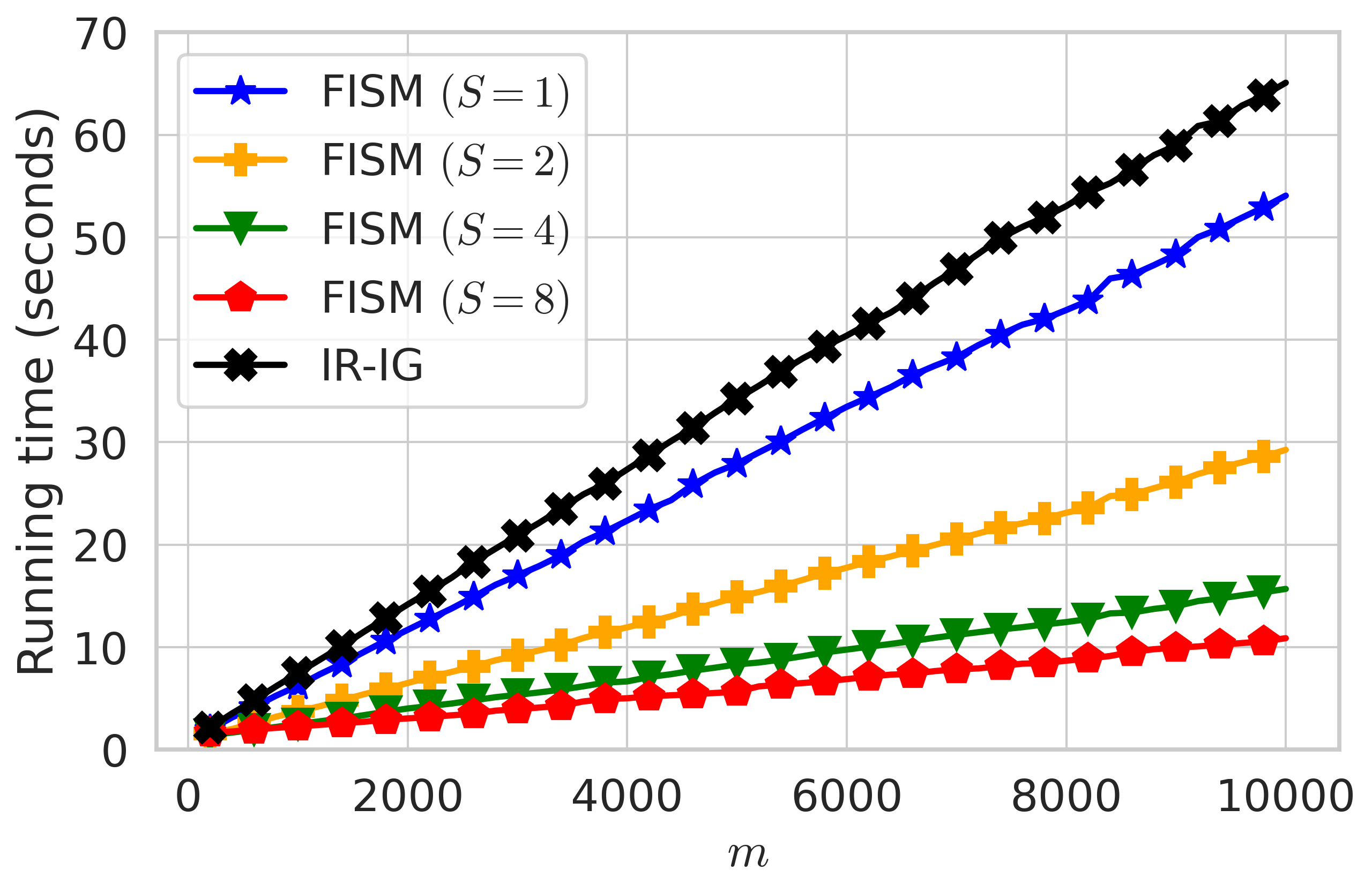}
    \caption{Comparison of total number of observations $m$ and computational running time for FISM and IR-IG.}
    \label{fig-m}
\end{figure}
\vspace{-0.2cm}

As demonstrated in Fig. \ref{fig-m},
FISM outperforms IR-IG for all considered values of $S$, with performance improving as $S$ increases.
This may be due to FISM's capacity for parallel computation and the ability of 
its groups
to process tasks independently, as discussed in Remark \ref{remark3.3}. Notably, even when $S=1$, FISM still achieve lower running times compared to IR-IG.
This is likely due to the updating procedure used in each iteration, as described in Remark \ref{remark3.1}.
\begin{figure}[H]
    \centering
    \begin{subfigure}[t]{0.23\textwidth}
        \centering
        \includegraphics[width=0.95\textwidth]{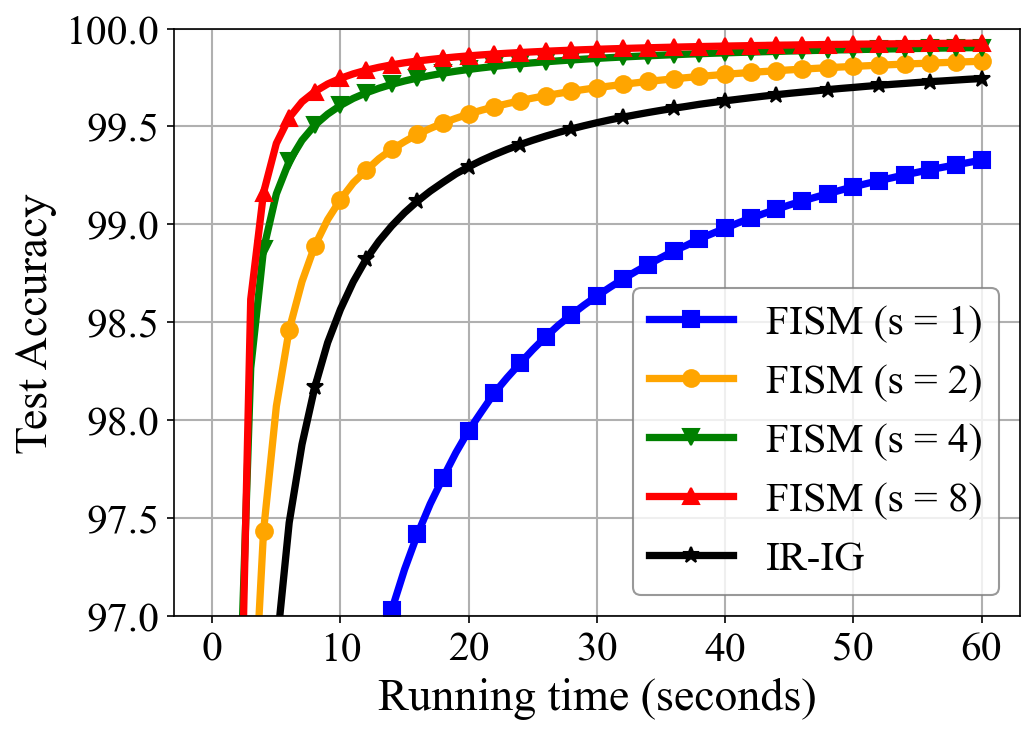}
        \caption{Test accuracy}
        \label{fig-acc}
    \end{subfigure}
    \hspace{-0.3cm}
    \begin{subfigure}[t]{0.23\textwidth}
        \centering
        \includegraphics[width=0.95\textwidth]{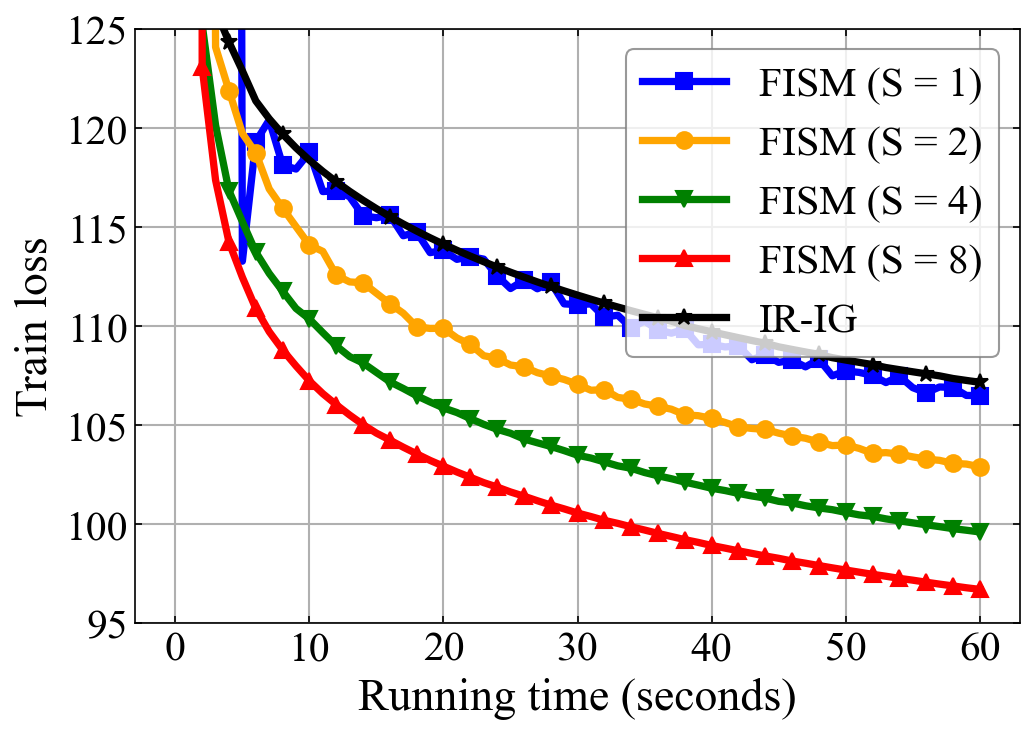}
        \caption{Train loss}
        \label{fig-loss}
    \end{subfigure}
    \caption{Comparison of performance for FISM and IR-IG}
    \label{fig:3}
\end{figure}
\vspace{-0.35cm}

In Fig. \ref{fig:3}, we present the test accuracy and train loss of FISM across different numbers of groups, and of IR-IG within $60$ seconds. 
We train the model using $m = 11000$ total training images and $2115$ testing images.
As shown,
 FISM with $S =4$ and $S=8$
 achieve higher test accuracy and lower train loss within shorter running times compared to FISM with $S=1$ and $S=2$, as well as to IR-IG. These results indicate that, in this experiment, increasing $S$ leads to better performance. However, Fig. \ref{fig:3} reveals that the performance of FISM with $S = 1$ is still worse than IR-IG, although it is faster in terms of running time. We suspect this stems from how we defined the term $\mathcal{H}_k$, which contrasts with IR-IG's approach of using the most recent parameter to compute the subgradient of $H$.
\vspace{-0.3cm}
\subsection{Location Problem}
\vspace{-0.1cm}
In this experiment, we focus on the location problem with a given anchor point and a finite number of target sets: 

\vspace{-0.3cm}
{\small \begin{equation*}
\begin{aligned}
    &\textrm{minimize} \quad && H(x):= 0.5 \| x - a \|^2 \\
    &\textrm{subject to} \quad && x \in {\textrm{argmin}}_{y \in X}  {\small F(y):=\sum_{i=1} ^S \sum_{j_i = 1}^{I_i} \textrm{dist}(y,C_{j_i}^i)},
\end{aligned}
\end{equation*}
}
\vspace{-0.2cm}

\noindent where the target sets $C_{j_i}^ i :=\{u\in \mathbb{R}^n:\|u-c_{j_i}^i\|\leq r_{j_i}^i\}$ are closed balls for all $i \in [S], j_i \in [I_i]$, with $m=\sum_{i=1}^S I_i$, and the anchor point $a$ is fixed. The main feature of this problem is to select the closest solution to the given anchor point $a$ that minimizes the sum of distances to all target sets $C_{j_i}^i$. 

Here, we investigate the performance of FISM with different numbers of groups ($S = 1, 2, 4,8$) and compare it with that of IR-IG under the following setting: the constraint set $X:=[-10, 10]^n$ is the box; the center point $c_{j_i}^i$ and the anchor point $a$, along with the initial point, are independently sampled in $(-10, 10)^n$; the radius $r_{j_i}^i$ are sampled between $(0, 1)$. We set the stopping criterion as $\textrm{max} \left\{ \frac{\| x_{k+1} - x_k\|}{\| x_k \|+1} , \frac{\left\vert H(x_{k+1}) - H(x_k) \right\vert}{H(x_k) +1},\frac{\left\vert F(x_{k+1}) - F(x_k) \right\vert}{F(x_k) +1}  \right\} \leq 10^{-5}$. 
In a similar fashion as previous experiment, we choose $ \gamma_k := \frac{1}{k^{0.8}}$  and $\lambda_k := \frac{1}{k^{0.1}}$. We report the average computational running times (in seconds) over 10 independent tests
across various problem dimensions ($n$) and numbers of target sets ($m$).
The results are shown in Table \ref{tab1}. 

\begin{table}[H]
\caption{
Average computational running times (in seconds) of IR-IG and FISM for $S=1,2,4$ and $8$. 
}
\label{tab1}
\centering
{ \scriptsize 

\begin{tabular}{l l r r r r r r r r r r}
\toprule
$n$
& $m$
& IR-IG 
& $S=1$
& $S=2$
& $S=4$
& $S=8$ \\
\midrule
10
 & 500   & 16.92   & 15.96   & 5.95   & 2.34  & 1.39 \\
 & 1000  & 42.88   & 40.51   & 13.47  & 5.31   & 2.76 \\
 & 5000  & 374.42  & 351.27  & 120.96  & 42.67   & 18.12 \\
 & 10000 & 849.49  & 800.44  & 295.64  & 104.30  & 38.62 \\
\midrule
50
 & 500   & 15.78   & 14.77  & 6.13   & 3.37   & 2,77 \\
 & 1000  & 35.47   & 33.20  & 12.46   & 5.06   & 2.96\\
 & 5000  & 318.56  & 298.60  & 106.08 & 36.96  & 15.26 \\
 & 10000 & 808.79  & 760.02  & 277.83  & 98.70  & 36.61 \\
\midrule
100
 & 500   & 17.33   & 16.17  & 8.31  &  4.74  &  3.91\\
 & 1000  & 34.73  & 32.44  &  12.90 &  5.95  & 4.10\\
 & 5000  & 292.83  & 276.16  & 97.92 &  35.52  & 14.08 \\
 & 10000 & 768.78  & 724.63 &  256.56 &  89.95 &35.93  \\
\bottomrule
\end{tabular}}
\end{table}

\vspace{-0.1cm}

In Table \ref{tab1}, FISM outperforms IR-IG across all values of $S$. 
Increasing $S$ lead to better performance, with $S=4$ and $S=8$ yielding the best results. This trend becomes even more apparent as the problem size increases. Notably, when
$m=10000$, FISM with $S=8$ reaches the stopping criteria at least $20$ times faster than IR-IG. 

\vspace{-0.2cm}
\section{Concluding Remarks}

In this letter, we presented an iterative method for BOPs within a FL framework. 
The method preserves the privacy and security of each client’s data, as well as the central server’s data, while providing computational savings.
It can be observed from Remark \ref{remark3.3} that in ideal situations, where all clients require the same computational time and communication costs are negligible, the time required per iteration for FISM could be up to $S$ times less than that of IR-IG. However, $T_{\text{FISM}}$ also depends on $\textrm{max}_{i \in [S]}\sum_{j_i = 1}^{I_i} s_{i,j_i}$, which introduces a bottleneck that can cause the system to wait for the slowest client.
 In heterogeneous scenarios, this issue consequently increases the time required for each iteration. To mitigate this issue, incorporating
asynchronous communication strategies would be an interesting direction for future work.

Furthermore, the example discussed in Remark \ref{remarkss} requires the tuning of four hyperparameters, leading to significant computational overhead. This comes from the fact that each iteration demands a complete pass through all local data points.
This computational challenge also motivates the exploration of adaptive stochastic variances in future work.

Finally, it is worth noting that objective functions in many practical applications frequently violate convexity assumptions in both the lower-level and upper-level problems, extending this framework to nonconvex settings represents a direction for enhancing its applicability to real-world scenarios.

\vspace{-0.2cm}

\section*{Acknowledgment}
The authors would like to thank the editor and the three anonymous referees for their valuable comments and remarks, which improved the quality and presentation of this work. The authors also thank Tipsuda Arunrat for her helpful suggestions.

\vspace{-0.3cm}

\bibliographystyle{IEEEtran}
\bibliography{references2}
\end{document}